\documentclass[A4j,11pt]{article}
\usepackage{tikz}

\usepackage{amssymb,amsmath,amsthm,amscd,amsfonts}
\usepackage{amsmath,amssymb,amsthm,amscd,ascmac,mathrsfs,bm,type1cm,multirow,array,enumerate,mathrsfs,tabularx}
\usepackage{graphicx}
\usepackage{graphics}

\newcommand{\ve}{\bm{e}} 

\newcommand{\vv}{{\bm{v}}} 
\newcommand{\vw}{{\bm{w}}} 

\begin{document}

\title{{\bf Isoparametric submanifolds\\
in Hilbert spaces\\
and holonomy maps}}
\author{{\bf Naoyuki Koike}}
\date{}
\maketitle

\abstract{
Let $\pi:P\to B$ be a smooth $G$-bundle over a compact Riemannian manifold $B$ and $c$ a smooth loop in $B$ 
of constant seed $a(>0)$, where $G$ is compact semi-simple Lie group.  In this paper, we prove that the holonomy map ${\rm hol}_c:\mathcal A_P^{H^s}\to G$ 
is a homothetic submersion of coefficient $a$, where $s$ is a non-negative integer, $\mathcal A_P^{H^s}$ is the Hilbert space of all $H^s$-connections of the bundle $P$.  
In particular, we prove that, if $s=0$, then ${\rm hol}_c$ has minimal regularizable fibres.  
From this fact, we can derive that each component of the inverse image of any equifocal submanifold in $G$ 
by the holonomy map ${\rm hol}_c:\mathcal A_P^{H^0}\to G$ is an isoparametric submanifold in $\mathcal A_P^{H^0}$.  
As the result, we obtain a new systematic construction of isoparametric submanifolds in a Hilbert space.  }

\section{Introduction}
R. S. Palais and C. L. Terng (\cite{PaTe}) introduced the notion of a {\it proper Fredholm submanifold} in a (separable) Hilbert space 
as an immersed (Hilbert) submanifold of finite codimension satisfying the following conditions:

\vspace{0.15truecm}

(PF-i)\ The normal exponential map $\exp^{\perp}$ of the submanifold is a Fredholm map;

(PF-ii)\ The restriction of $\exp^{\perp}$ to the ball normal bundle of any radius is proper.  

\vspace{0.15truecm}

\noindent
Since the shape operators of a proper Fredholm submanifold are self-adjoint and compact operators, the set of all eigenvalues of each shape operator 
is bounded and has no accumulating point other than 0.  Also, the multiplicities of the nonzero eigenvalues are finite.  
Hence, for each unit normal vector $\xi$ of the submanifold, the set of of all focal radii (which are equal to the inverses of the nonzero eignenvalues of 
the shape operator $A_{\xi}$) of the submanifold along the normal geodesic $\gamma_{\xi}$ of $\xi$-direction has no accumulating points and the multiplicities 
of each focal radii are finite.  Thus a proper Fredholm submanifold has a good focal structure similar to finite dimensional submanifolds.  
Furthermore, they (\cite{PaTe}) defined the notion of an {\it isoparametric submanifold} in a Hilbert space as a proper Fredholm submanifold satisfying the following conditions:

\vspace{0.15truecm}

(I-i)\ The normal holonomy group of the submanifold is trivial;

(I-ii)\ For any parallel normal vector field $\widetilde{\xi}$ of the submanifold, the set of all eigenvalues of the shape operator $A_{\widetilde{\xi}_p}$ 
is independent of the choice of the base point $p$ with considering the multiplicities.  

\vspace{0.15truecm}

\noindent
See \cite{GH}, \cite{H}, \cite{HL}, \cite{HLO}, \cite{PaTe}, \cite{PiTh}, \cite{T1},\cite{T2} and \cite{TT} 
about the study of isoparametric submanifolds in a Hilbert space.  

Next we rcall the definitions of the regularized traces of a self-adjoint and compact operator of a (separable) Hilbert space in the sense of \cite{KT} and \cite{HLO}.  
Let $A$ be a self-adjoint operator of a (separable) Hilbert space $(V,\langle\,\,,\,\,\rangle)$ 
and 
$$-\lambda^-_1\leq-\lambda^-_2\leq\cdots\leq 0\leq\cdots\leq\lambda^+_2\leq\lambda^+_1$$
be the spectrum of $A$.  
C. King and C. L. Terng (\cite{KT}) defined the $L^s$-norm $\|A\|_s$ ($s>1$) of of a self-adjoint and compact operator $A$ by 
$$\|A\|_s:=\left(\sum\limits_{i=1}^{\infty}((\lambda^+_i)^s+(\lambda^-_i)^s)\right)^{\frac{1}{s}}.$$
They called that $A$ is {\it regularizable} if $\|A\|_s<\infty$ for all $s>1$ and if 
$$\lim_{s\downarrow 1}\sum\limits_{i=1}^{\infty}((\lambda^+_i)^s-(\lambda^-_i)^s)\leqno{(1.1)}$$
exists.  Also, they called the limit in $(1.1)$ the {\it regularized trace} of $A$.  
In \cite{HLO}, this regularized trace was called the {\it $\zeta$-regularized trace}.  
Later, E. Heintze, X. Liu and C. Olmos (\cite{HLO}) called that $A$ is {\it regulaizable} if 
$\|A\|_2<\infty$ and 
$$\sum\limits_{i=1}^{\infty}(\lambda^+_i-\lambda^-_i)\leqno{(1.2)}$$
exists.  Also, they called the limit in $(1.2)$ the {\it regularized trace} of $A$.  
The regularized trace in the sense of \cite{HLO} is easier to handle than one in 
the sense of \cite{KT}.  
We denote the regularized trace (in the sense of \cite{HLO}) (resp. the $\zeta$-regularized trace) of $A$ by ${\rm Tr}_r\,A$ (resp. ${\rm Tr}_{\zeta}\,A$), 
where the subscript ``$r$'' in ${\rm Tr}_r$ implies the initial letter of ``regularized''.  
In almost all relevant cases, these regularized traces coincide.  
In this paper, we shall use the terminologies ''regularized trace'' and ``regularizable'' in the sense of \cite{HLO}.  
Here we give an example of a self-ajoint and compact operator $A$ such that $\|A\|_s$ ($s>1$) and ${\rm Tr}_{\zeta}\,A$ exist but ${\rm Tr}_r\,A$ does not exist.  

\vspace{0.25truecm}

\noindent
{\it Example}\ \ We consider the case where the spectrum of $A$ is given by 
$$-1<-\frac{1}{2}<\cdots<-\frac{1}{k}<\cdots\leq 0\leq\cdots<\frac{2}{k}<\cdots<1<2,$$
where we note that the multiplcities of all eigenvalues of $A$ are equal to $1$.  
Then we have 
\begin{align*}
\|A\|_s=&\sum_{k=1}^{\infty}\left(\left(\frac{2}{k}\right)^s+\left(\frac{1}{k}\right)^s\right)=\sum_{k=1}^{\infty}\frac{2^s+1}{k^s}\\
=&(2^s+1)\left(1+\sum_{k=2}^{\infty}\frac{1}{k^s}\right)<(2^s+1)\left(1+\int_1^{\infty}\frac{1}{x^s}\,dx\right)\\
=&\frac{s(2^s+1)}{s-1}\qquad\qquad\qquad(s>1)
\end{align*}
and 
\begin{align*}
{\rm Tr}_{\zeta}\,A=&\lim_{s\downarrow 1}\sum_{k=1}^{\infty}\left(\left(\frac{2}{k}\right)^s-\left(\frac{1}{k}\right)^s\right)
=\lim_{s\downarrow 1}\sum_{k=1}^{\infty}\frac{2^s-1}{k^s}\\
=&\lim_{s\downarrow 1}(2^s-1)\left(1+\sum_{k=2}^{\infty}\frac{1}{k^s}\right)
<\lim_{s\downarrow 1}(2^s-1)\left(1+\int_1^{\infty}\frac{1}{x^s}\,dx\right)\\
=&\lim_{s\downarrow 1}\frac{s(2^s-1)}{s-1}=1+2\ln\,2.
\end{align*}
Thus both $\|A\|_s$ ($s>1$) and ${\rm Tr}_{\zeta}\,A$ exists.  Hence $A$ is regularizable in the sense of \cite{KT}.  
On the other hand, we have 
$${\rm Tr}_r\,A=\sum_{k=1}^{\infty}\left(\frac{2}{k}-\frac{1}{k}\right)=\sum_{k=1}^{\infty}\frac{1}{k}=\infty.$$
Thus ${\rm Tr}_r\,A$ does not exist.  Hence $A$ is not regularizable in the sense of \cite{HLO}.  

\vspace{0.25truecm}

Let $M$ be a proper Fredholm submanifold in $(V,\langle\,\,,\,\,\rangle)$ immersed by $f$.  
If, for any $p\in M$ and any normal vector $\xi$ of $M$ at $p$, there exist the regularized trace ${\rm Tr}_r\,(A_p)_{\xi}$ of the shape operator 
$(A_p)_{\xi}$ of $M$ for $\xi$ and the trace ${\rm Tr}\,(A_p)_{\xi}^2$ of $(A_p)_{\xi}^2$, then $M$ is called a {\it regularizable submanifold}.  
In particular, if, for any $p\in M$ and any normal vector $\xi$ of $M$ at $p$, ${\rm Tr}_r\,(A_p)_{\xi}$ vanishes, then $M$ is called a 
{\it minimal regularizable submanifold}.  

Y. Maeda, S. Rosenberg and P. Tondeur (\cite{MRT}) studied the submanifold structure of the gauge orbit in the Hilbert space of 
all $H^0$-connections of a $G$-bundle $P$ over a compact Riemannian manifold $(B,g_B)$, where $G$ is a compact semi-simple Lie group.  
Denote by $\mathcal A_P^{H^0}$ the Hilbet space of all $H^0$-connections of $P$ and $\mathcal G_P^{H^1}$ the $H^1$-gauge transformation group of $P$.  
Note that, if ${\rm dim}\,B\geq 2$, then $\mathcal G_P^{H^1}$ is not a smooth Hilbert Lie group and 
the gauge action $\mathcal G_P^{H^1}\curvearrowright\mathcal A_P^{H^0}$ is not smooth.  
Take a $C^{\infty}$-curve $c:[0,1]\to B$ and a $C^{\infty}$-curve $\sigma:[0,1]\to P$ with $\pi\circ\sigma=c$.  
Let $c^{\ast}P$ be the pull-back bundle of $P$ by $c$, which is identified with the trivial $G$-bundle $[0,1]\times G$ over $[0,1]$ by $\sigma$.  
Hence the Hilbert space of all $H^0$-connections of $c^{\ast}P$ is identified with the Hilbert space $H^0([0,1],\mathfrak g)$ of all $H^0$-curves 
in the Lie algebra $\mathfrak g$ of $G$ by $\sigma$.  
Also, the $H^1$-gauge transformation group of $c^{\ast}P$ is identified with the Hilbert Lie group $H^1([0,1],G)$ of all $H^1$-curves in $G$.  
The {\it parallel transport map} for $G$ is a map $\phi:H^0([0,1],\mathfrak g)\to G$ defined by 
$$\phi(u):={\bf g}_u(1)\quad\,\,(u\in H^0([0,1],\mathfrak g)),$$
where ${\bf g}_u$ is the element of $H^1([0,1],G)$ with ${\bf g}_u(0)=e$ ($e\,:\,$ the identity element of $G$) and 
$(R_{{\bf g}_u(t)})_{\ast}^{-1}({\bf g}_u'(t))=u(t)$ $(t\in[0,1])$.  
The group $H^1([0,1],G)$ acts on $H^0([0,1],\mathfrak g)$ as the (usual) action of the gauge transformation group on the space of the connections 
and the closed subgroup $\Lambda_e^{H^1}(G)$ (of $H^1([0,1],G)$) of all $H^1$-loops at $e$ in $G$ acts on $H^0([0,1],\mathfrak g)$ as its subaction.  
This subaction $\Lambda_e^{H^1}(G)\curvearrowright H^0([0,1],\mathfrak g)$ is transitive and the orbit space $H^0([0,1],\mathfrak g)/\Lambda_e^{H^1}(G)$ 
is identified with $G$.  Furthermore the orbit map of this subaction $\Lambda_e^{H^1}(G)\curvearrowright H^0([0,1],\mathfrak g)$ is identified with 
the parallel transport map $\phi$.  
Fix an ${\rm Ad}(G)$-invariant inner product $\langle\,\,,\,\,\rangle_{\mathfrak g}$ of $\mathfrak g$.  
The parallel transport map $\phi$ is a Riemannian submersion of the Hilbert space $H^0([0,1],\mathfrak g)$ equipped with the $H^0$-inner product 
defined by $\langle\,\,,\,\,\rangle_{\mathfrak g}$ onto $G$ equipped with a bi-invariant metric defined by $\langle\,\,,\,\,\rangle_{\mathfrak g}$.  
See \cite{M1}, \cite{M2}, \cite{PaTe}, \cite{PiTh}, \cite{T1}, \cite{T2}, \cite{TT} and \cite{K} about the study related to the parallel transport map.  

In a finite dimensional complete Riemannian manifold $(\widetilde M,\widetilde g)$, the notion of an {\it equifocal submanifold} is defined as an immersed compact submanifold 
satisfying the following conditions:

\vspace{0.15truecm}

(i)\ The submanifold has abelian normal bundle (i.e., flat section);

(ii)\ The normal holonomy group of the submanifold is trivial;

(iii)\ The focal structure of the submanifold is invariant under parallel translations with respect to the normal connection.  

\vspace{0.15truecm}

\noindent
According to Theorem 1.10 of \cite{TT}, if $M$ is an equifocal submanifold in $G$ equipped with the bi-invariant metric induced from a fixed ${\rm Ad}(G)$-invariant 
inner product $\langle\,\,,\,\,\rangle$ of $\mathfrak g$, then each component of $\phi^{-1}(M)$ is an isoparametric submanifold in $H^0([0,1],\mathfrak g)$ equipped with 
the $L^2$-inner product induced from $\langle\,\,,\,\,\rangle$.  

Assume that $\displaystyle{s>\frac{1}{2}\,{\rm dim}\,B-1}$.  
Let $\mathcal A_P^{H^s}$ be the Hilbert space of all $H^s$-connections of $P$ and $\mathcal G_P^{H^{s+1}}$ the group of all $H^{s+1}$-gauge transformations 
of $P$.  
Then, according to Section 9 of \cite{P}, $\mathcal G_P^{H^{s+1}}$ is a smooth Hilbert Lie group.  
Also, according to Lemma 1.2 of \cite{U}, the gauge action $\mathcal G_P^{H^{s+1}}\curvearrowright\mathcal A_P^{H^s}$ is smooth.  
However, $\mathcal G_P^{H^{s+1}}$ does not act isometrically on the Hilbert space $(\mathcal A_P^{H^s},\langle\,\,,\,\,\rangle^{\omega_0}_s)$, 
where $\langle\,\,,\,\,\rangle_s^{\omega_0}$ denotes the $L^2_s$-inner product defined by using an arbitrarily fixed $C^{\infty}$-connection $\omega_0$ of $P$ 
(see Section 2 about this defintition).  Define a Riemannian metric ${\it g}_s$ on $\mathcal A_P^{H^s}$ by 
$({\it g}_s)_{\omega}:=\langle\,\,,\,\,\,\rangle^{\omega}_s$ ($\omega\in\mathcal A_P^{H^s}$), where $\langle\,\,,\,\,\rangle_s^{\omega}$ denotes 
the $L^2_s$-inner product defined by using $\omega$ instead of $\omega_0$.  
This Riemannian metric ${\it g}_s$ is non-flat.  
The gauge transformation group $\mathcal G_P^{H^{s+1}}$ acts isometrically on the Riemannian Hilbert manifold $(\mathcal A_P^{H^s},{\it g}_s)$, 
where we note that the Hilbert space $(\mathcal A_P^{H^s},\langle\,\,,\,\,\rangle^{\omega_0}_s)$ is regarded as the tangent space of $(\mathcal A_P^{H^s},{\it g}_s)$ 
at $\omega_0$.  Hence we can give the moduli space $\mathcal M_P^{H^s}:=\mathcal A_P^{H^s}/\mathcal G_P^{s+1}$ the Riemannian orbimetric 
$\overline{\it g}_s$ such that the orbit map $\pi_{\mathcal M_P}:(\mathcal A_P^{H^s},{\it g}_s)\to(\mathcal M_P^{H^s},\overline{\it g}_s)$ is a Riemannian orbisubmersion.  

Fix an ${\rm Ad}(G)$-invariant inner product $\langle\,\,,\,\,\rangle_{\mathfrak g}$ of $\mathfrak g$.  
Let $\langle\,\,,\,\,\rangle_0$ be the $L^2$-inner product of $\mathcal A_P^{H^0}$ defined by using the Riemannian metric $g_B$ of $B$ and 
$\langle\,\,,\,\,\rangle_{\mathfrak g}$ and $g_G$ the bi-invariant metric of $G$ defined by using $\langle\,\,,\,\,\rangle_{\mathfrak g}$.  

In this paper, we prove the following fact.  

\vspace{0.25truecm}

\noindent
{\bf Theorem A.} {\sl Let $\pi:P\to B$ be a smooth $G$-bundle over a compact Riemannian manifold $(B,g_B)$ and 
$(\mathcal A_P^{H^s},g_s)$ the Riemannian Hilbert manifold consisting of all $H^s$-connections of the bundle $P$, where $s$ is an arbitrary non-negative number.  
Then, for any $C^{\infty}$-loop $c:S^1\to B$ of constant speed $a(>0)$, the following statements (i) and (ii) hold:

{\rm (i)}\  The holonomy map ${\rm hol}_c:(\mathcal A_P^{H^s},g_s)\to(G,g_G)$ along $c$ is a homothetic submersion of coefficient $a$.  
In particular, ${\rm hol}_c:(\mathcal A_P^{H^0},g_0)\to(G,g_G)$ has minimal regularizable fibres;

{\rm (ii)}\ For an equifocal submanifold $M$ in $(G,g_G)$, each component of \newline
${\rm hol}_c^{-1}(M)$ is a $(\mathcal G_P^{H^1})_{c(0)}$-invariant isoparametric submanifold in 
$(\mathcal A_P^{H^0},g_0)$, \newline
where $(\mathcal G_P^{H^1})_{c(0)}$ is the based gauge transformation group of class $H^1$ at $c(0)$ of the bundle $P$.  
Furthermore, if $M$ is ${\rm Ad}(G)$-invariant, then each component of ${\rm hol}_c^{-1}(M)$ is $\mathcal G_P^{H^1}$-invariant, where 
${\rm Ad}$ is the adjoint action of $G$ on oneself and $\mathcal G_P^{H^1}$ is the gauge transformation group of class $H^1$ of the bundle $P$.}

\section{The holonomy map and the pull-back connection map along a loop} 
In this section, we define the holonomy map along loops and the pull-back connection maps by loops for smooth principal bundle.  
Let $\pi:P\to B$ be a $G$-bundle of class $C^{\infty}$ over a compact Riemannian manifold $(B,g_B)$, where $G$ is a compact semi-simple Lie group.  
Fix an ${\rm Ad}(G)$-invariant inner product $\langle\,\,,\,\,\rangle_{\mathfrak g}$ (for example, the $(-1)$-multiple of the Killing form of $\mathfrak g$) of 
the Lie algebra $\mathfrak g$ of $G$, where ${\rm Ad}$ denotes the adjoint representation of $G$.  Denote by $g_G$ the bi-invariant metric of $G$ induced from 
$\langle\,\,,\,\,\rangle_{\mathfrak g}$.  Let $\mathcal A_P^{\infty}$ be the affine Hilbert space of all $C^{\infty}$-connections of $P$ 
and $\Omega^{\infty}_{\mathcal T,i}(P,\mathfrak g)$ be the Hilbert space of all $\mathfrak g$-valued tensorial $i$-form of class $C^{\infty}$ on $P$ 
(see \cite{KN} about the definition of a $\mathfrak g$-valued tensorial $i$-form).
Also, let 
$$\Omega_i^{\infty}(B,{\rm Ad}(P))(=\Gamma^{\infty}((\wedge^iT^{\ast}B)\otimes{\rm Ad}(P)))$$
be the space of all ${\rm Ad}(P)$-valued $i$-forms of class $C^{\infty}$ over $B$, where ${\rm Ad}(P)$ denotes the adjoint bundle 
$P\times_{\rm Ad}\mathfrak g$.  The space ${\mathcal A}_P^{\infty}$ is the affine space having $\Omega^{\infty}_{\mathcal T,1}(P,\mathfrak g)$ as the associated vector space.  
Furthermore, $\Omega^{\infty}_{\mathcal T,1}(P,\mathfrak g)$ is identified with $\Omega_1^{\infty}(B,{\rm Ad}(P))$ under the correspondence $A\leftrightarrow\widehat A$ 
defined by $u\cdot A_u(\vv)=\widehat A_{\pi(u)}(\pi_{\ast}(\vv))$ ($u\in P,\,\,\,\vv\in T_uP$), 
where we note that $\widehat A$ is the section of ${\rm Ad}(P)$ defined from the pullbacks of $A$ by local sections of $P$.  

According to \cite{GP2}, we shall define the $H^s$-completion of $\Omega_{\mathcal T,1}^{\infty}(P,\mathfrak g)$, where $s\geq 0$.  
Denote by $\mathcal A_P^{w,s}$ the space of all $s$-times weak differentiable connections of $P$ and 
$\Omega_{\mathcal T,i}^{w,s}(P,\mathfrak g)$ the space of all $s$-times weak differentiable tensorial $i$-forms on $P$.  
Fix a $C^{\infty}$-connection $\omega_0$ of $P$ as the base point of $\mathcal A_P^{w,s}$.  
Define an operator $\square_{\omega_0}:\Omega_{\mathcal T,i}^{w,s}(P,\mathfrak g)\to\Omega_{\mathcal T,i}^{w,s-2}(P,\mathfrak g)$ by 
$$
\square_{\omega_0}:=\left\{\begin{array}{ll}
d_{\omega_0}\circ d_{\omega_0}^{\ast}+d_{\omega_0}^{\ast}\circ d_{\omega_0}+{\rm id} & (i\geq 1)\\
d_{\omega_0}^{\ast}\circ d_{\omega_0}+{\rm id} & (i=0),
\end{array}\right.$$
where $d_{\omega_0}$ denotes the covariant exterior derivative with respect to $\omega_0$ and 
$d_{\omega_0}^{\ast}$ denotes the adjoint operator of $d_{\omega_0}$ with respect to the $L^2$-inner products of 
$\Omega_i^{w,j}(B,{\rm Ad}(P))$ ($i\geq 0,\,\,\,j\geq 1$).  
The $L^2_s$-inner product $\langle\,\,,\,\,\rangle^{\omega_0}_s$ of 
$$T_{\omega}\mathcal A_P^{w,s}(\approx\Omega_{\mathcal T,1}^{w,s}(P,\mathfrak g)\approx\Omega_1^{w,s}(B,{\rm Ad}(P))\approx\Gamma^{w,s}(T^{\ast}B\otimes{\rm Ad}(P)))$$
is defined by 
$$\begin{array}{r}
\displaystyle{\langle A_1,A_2\rangle^{\omega_0}_s
:=\int_{x\in B}\langle(\widehat A_1)_x,(\widehat{\square_{\omega_0}^s(A_2)})_x\rangle_{B,\mathfrak g}\,dv_B}\\
\displaystyle{(A_1,A_2\in\Omega_{\mathcal T,1}^{w,s}(P,\mathfrak g)),}
\end{array}\leqno{(2.1)}$$
where $\widehat{\square_{\omega_0}^s(A_2)}$ denotes the ${\rm Ad}(P)$-valued $1$-form over $B$ corresponding to 
$\square_{\omega_0}^s(A_2)$, $\langle\,\,,\,\,\rangle_{B,\mathfrak g}$ denotes the fibre metric of 
$T^{\ast}B\otimes{\rm Ad}(P)$ defined naturally from $g_B$ and $\langle\,\,,\,\,\rangle_{\mathfrak g}$ and $dv_B$ denotes the volume element of $g_B$.  
Let $\Omega_{\mathcal T,1}^{H^s}(P,\mathfrak g)$ be the completion of $\Omega_{\mathcal T,1}^{\infty}(P,\mathfrak g)$ with respect to 
$\langle\,\,,\,\,\rangle^{\omega_0}_s$.  
Also, set 
$${\mathcal A}_P^{H^s}:=\{\omega_0+A\,\vert\,A\in\Omega_{\mathcal T,1}^{H^s}(P,\mathfrak g)\}.$$
Let $\Omega_1^{H^s}(B,{\rm Ad}(P))$ be the completion of 
$\Omega_1^{\infty}(B,{\rm Ad}(P))$ corresponding to $\Omega_{\mathcal T,1}^{H^s}(P,\mathfrak g)$.  

\vspace{0.25truecm}

\centerline{
{\small
\unitlength 0.1in
\begin{picture}( 48.2000,  6.4000)( -9.0000,-15.6000)
\put(43.6000,-11.4000){\makebox(0,0)[rb]{$\Omega_1^{H^s}(B,{\rm Ad}(P))$}}%
\put(14.5000,-11.6000){\makebox(0,0)[lb]{$T_{\omega_0}\mathcal A_P^{H^s}=\Omega_{\mathcal T,1}^{H^s}(P,\mathfrak g)$}}%
\put(10.7000,-10.9000){\makebox(0,0)[lb]{$\approx$}}%
%
\special{pn 8}%
\special{ar 3860 1210 60 90  6.2831853 6.2831853}%
\special{ar 3860 1210 60 90  0.0000000 3.1415927}%
%
\special{pn 8}%
\special{pa 3860 1300}%
\special{pa 3860 1210}%
\special{fp}%
\put(38.0000,-13.9000){\makebox(0,0)[lt]{$\widehat A$}}%
\put(17.4000,-13.6000){\makebox(0,0)[lt]{$A(:=\omega-\omega_0)$}}%
%
\special{pn 8}%
\special{ar 1800 1210 60 90  6.2831853 6.2831853}%
\special{ar 1800 1210 60 90  0.0000000 3.1415927}%
%
\special{pn 8}%
\special{pa 1800 1300}%
\special{pa 1800 1210}%
\special{fp}%
\put(9.0000,-11.5000){\makebox(0,0)[rb]{$\mathcal A_P^{H^s}$}}%
\put(6.9000,-13.9000){\makebox(0,0)[lt]{$\omega$}}%
%
\special{pn 8}%
\special{ar 740 1210 60 90  6.2831853 6.2831853}%
\special{ar 740 1210 60 90  0.0000000 3.1415927}%
%
\special{pn 8}%
\special{pa 740 1300}%
\special{pa 740 1210}%
\special{fp}%
%
\special{pn 8}%
\special{pa 1292 1446}%
\special{pa 1040 1446}%
\special{fp}%
\special{sh 1}%
\special{pa 1040 1446}%
\special{pa 1108 1466}%
\special{pa 1094 1446}%
\special{pa 1108 1426}%
\special{pa 1040 1446}%
\special{fp}%
%
\special{pn 8}%
\special{pa 1300 1446}%
\special{pa 1550 1446}%
\special{fp}%
\special{sh 1}%
\special{pa 1550 1446}%
\special{pa 1484 1426}%
\special{pa 1498 1446}%
\special{pa 1484 1466}%
\special{pa 1550 1446}%
\special{fp}%
%
\special{pn 8}%
\special{pa 3360 1450}%
\special{pa 2784 1450}%
\special{fp}%
\special{sh 1}%
\special{pa 2784 1450}%
\special{pa 2852 1470}%
\special{pa 2838 1450}%
\special{pa 2852 1430}%
\special{pa 2784 1450}%
\special{fp}%
%
\special{pn 8}%
\special{pa 3340 1450}%
\special{pa 3530 1450}%
\special{fp}%
\special{sh 1}%
\special{pa 3530 1450}%
\special{pa 3464 1430}%
\special{pa 3478 1450}%
\special{pa 3464 1470}%
\special{pa 3530 1450}%
\special{fp}%
\put(29.9000,-10.9000){\makebox(0,0)[lb]{$\approx$}}%
\put(15.8000,-15.6000){\makebox(0,0)[lt]{($\omega_0:$ the base point of $\mathcal A_P^{H^0}$)}}%
\end{picture}%
}
\hspace{4truecm}}

\vspace{0.85truecm}

Let $\mathcal G_P^{\infty}$ be the group of all $C^{\infty}$-gauge transformations ${\bf g}$'s of $P$ with $\pi\circ{\bf g}=\pi$.  
For each ${\bf g}\in\mathcal G_P^{\infty}$, 
$\widehat{\bf g}\in C^{\infty}(P,G)$ is defined by ${\bf g}(u)=u\widehat{\bf g}(u)\,\,\,(u\in P)$.  
This element $\widehat{\bf g}$ satisfies 
$$\widehat{\bf g}(ug)={\rm Ad}(g^{-1})(\widehat{\bf g}(u))\,\,\,\,(\forall u\in P,\,\,\forall g\in G),$$
where ${\rm Ad}$ denotes the homomorphism of $G$ to ${\rm Aut}(G)$ defined by 
${\rm Ad}(g_1)(g_2):=g_1\cdot g_2\cdot g_1^{-1}$ ($g_1,g_2\in G$).  
Under the correspondence ${\bf g}\leftrightarrow\widehat{\bf g}$, $\mathcal G_P^{\infty}$ is identified with 
$$\widehat{\mathcal G}_P^{\infty}:=\{\widehat{\bf g}\in C^{\infty}(P,G)\,|\,
\widehat{\bf g}(ug)={\rm Ad}(g^{-1})(\widehat{\bf g}(u))\,\,\,(\forall\,u\in P,\,\,\forall\,g\in G)\}.$$
For $\widehat{\bf g}\in\widehat{\mathcal G}_P^{\infty}$, 
the $C^{\infty}$-section $\breve{\bf g}$ of the associated $G$-bundle $P\times_{{\rm Ad}}G$ is 
defined by 
$\breve{\bf g}(x):=u\cdot\widehat{\bf g}(u)\,\,\,(x\in B)$, where $u$ is any element of $\pi^{-1}(x)$.  
Under the correspondence $\widehat{\bf g}\leftrightarrow\breve{\bf g}$, 
$\widehat{\mathcal G}_P^{\infty}(=\mathcal G_P^{\infty})$ is identified with the space 
$\Gamma^{\infty}(P\times_{{\rm Ad}}G)$ of all $C^{\infty}$-sections of $P\times_{{\rm Ad}}G$.  
The gauge action $\mathcal G_P^{\infty}\curvearrowright\mathcal A_P^{\infty}$ is given by 
$$\begin{array}{r}
\displaystyle{({\bf g}\cdot\omega)_u={\rm Ad}(\widehat{\bf g}(u))\circ\omega_u-(R_{\widehat{\bf g}(u)})_{\ast}^{-1}\circ\widehat{\bf g}_{\ast u}}\\
\displaystyle{({\bf g}\in\mathcal G_P^{\infty},\,\,\,\omega\in\mathcal A_P^{\infty}).}
\end{array}\leqno{(2.2)}$$
The $H^{s+1}$-completion of $\Gamma^{\infty}(P\times_{{\rm Ad}}G)$ was defined in \cite{GP1} (see Section 1 (P668)).  
Denote by $\Gamma^{H^{s+1}}(P\times_{{\rm Ad}}G)$ this $H^{s+1}$-completion.  
Also, denote by $\mathcal G_P^{H^{s+1}}$ (resp. $\widehat{\mathcal G}_P^{H^{s+1}}$) the $H^{s+1}$-completion of 
$\mathcal G_P^{\infty}$ (resp. $\widehat{\mathcal G}_P^{\infty}$) corresponding to 
$\Gamma^{H^{s+1}}(P\times_{{\rm Ad}}G)$.  
If $\displaystyle{s>\frac{1}{2}\,{\rm dim}\,B-1}$, then the $H^{s+1}$-gauge transformation group $\mathcal G_P^{H^{s+1}}$ of $P$ is a smooth Hilbert Lie group 
and the gauge action $\mathcal G_P^{H^{s+1}}\curvearrowright\mathcal A_P^{H^s}$ is smooth as stated in Introduction.  
However, by this action, $\mathcal G_P^{H^{s+1}}$ does not act isometrically on the Hilbert space 
$(\mathcal A_P^{H^s},\langle\,\,,\,\,\rangle^{\omega_0}_s)$.  
Define a Riemannian metric ${\it g}_s$ on $\mathcal A_P^{H^s}$ by 
$({\it g}_s)_{\omega}:=\langle\,\,\,\,\,\rangle^{\omega}_s$ ($\omega\in\mathcal A_P^{H^s}$), where 
$\langle\,\,,\,\,\rangle_s^{\omega}$ is the $L^2_s$-inner product defined as in $(2.1)$ by using $\omega$ instead of $\omega_0$.  
This Riemannian metric ${\it g}_s$ is non-flat.  
The gauge transformation group $\mathcal G_P^{H^{s+1}}$ acts isometrically on 
the Riemannian Hilbert manifold $(\mathcal A_P^{H^s},{\it g}_s)$, where we note that 
the Hilbert space $(\mathcal A_P^{H^s},\langle\,\,,\,\,\rangle^{\omega_0}_s)$ is regarded as 
the tangent space of $(\mathcal A_P^{H^s},{\it g}_s)$ at $\omega_0$.  
Hence we can give the moduli space $\mathcal M_P^{H^s}:=\mathcal A_P^{H^s}/\mathcal G_P^{s+1}$ the Riemannian orbimetric 
$\overline{\it g}_s$ such that the orbit map 
$\pi_{\mathcal M_P}:(\mathcal A_P^{H^s},{\it g}_s)\to(\mathcal M_P^{H^s},\overline{\it g}_s)$ is a Riemannian orbisubmersion.  

\vspace{0.25truecm}

\centerline{
{\small
\unitlength 0.1in
\begin{picture}( 23.6000,  4.8500)( 10.3000,-16.1500)
\put(10.3000,-13.2000){\makebox(0,0)[lb]{$\mathcal G^{H^{s+1}}_P$}}%
\put(15.4000,-13.1000){\makebox(0,0)[lb]{$\approx$}}%
\put(20.2000,-13.3000){\makebox(0,0)[lb]{$\widehat{\mathcal G}^{H^{s+1}}_P$}}%
\put(25.5000,-13.0000){\makebox(0,0)[lb]{$\approx$}}%
\put(29.7000,-13.4000){\makebox(0,0)[lb]{$\Gamma^{H^{s+1}}(P\times_{\rm Ad}G)$}}%
\put(11.1000,-15.3000){\makebox(0,0)[lt]{${\bf g}$}}%
\put(20.7000,-15.2000){\makebox(0,0)[lt]{$\widehat{\bf g}$}}%
\put(32.9000,-15.2000){\makebox(0,0)[lt]{$\breve{\bf g}$}}%
%
\special{pn 8}%
\special{ar 1150 1360 60 90  6.2831853 6.2831853}%
\special{ar 1150 1360 60 90  0.0000000 3.1415927}%
%
\special{pn 8}%
\special{pa 1150 1450}%
\special{pa 1150 1360}%
\special{fp}%
%
\special{pn 8}%
\special{ar 2110 1350 60 90  6.2831853 6.2831853}%
\special{ar 2110 1350 60 90  0.0000000 3.1415927}%
%
\special{pn 8}%
\special{pa 2110 1440}%
\special{pa 2110 1350}%
\special{fp}%
%
\special{pn 8}%
\special{ar 3330 1360 60 90  6.2831853 6.2831853}%
\special{ar 3330 1360 60 90  0.0000000 3.1415927}%
%
\special{pn 8}%
\special{pa 3330 1450}%
\special{pa 3330 1360}%
\special{fp}%
%
\special{pn 8}%
\special{pa 1692 1606}%
\special{pa 1430 1606}%
\special{fp}%
\special{sh 1}%
\special{pa 1430 1606}%
\special{pa 1498 1626}%
\special{pa 1484 1606}%
\special{pa 1498 1586}%
\special{pa 1430 1606}%
\special{fp}%
%
\special{pn 8}%
\special{pa 1700 1606}%
\special{pa 1920 1606}%
\special{fp}%
\special{sh 1}%
\special{pa 1920 1606}%
\special{pa 1854 1586}%
\special{pa 1868 1606}%
\special{pa 1854 1626}%
\special{pa 1920 1606}%
\special{fp}%
%
\special{pn 8}%
\special{pa 2810 1610}%
\special{pa 2476 1610}%
\special{fp}%
\special{sh 1}%
\special{pa 2476 1610}%
\special{pa 2542 1630}%
\special{pa 2528 1610}%
\special{pa 2542 1590}%
\special{pa 2476 1610}%
\special{fp}%
%
\special{pn 8}%
\special{pa 2820 1610}%
\special{pa 3030 1610}%
\special{fp}%
\special{sh 1}%
\special{pa 3030 1610}%
\special{pa 2964 1590}%
\special{pa 2978 1610}%
\special{pa 2964 1630}%
\special{pa 3030 1610}%
\special{fp}%
\end{picture}%
}
\hspace{0.75truecm}}

\vspace{0.25truecm}

Let $\displaystyle{S^1:=\{e^{2\pi\sqrt{-1}t}\,|\,t\in[0,1]\}}$.  
Define $z:[0,1]\to S^1$ by $z(t):=e^{2\pi\sqrt{-1}t}$ ($t\in[0,1]$).  
Fix $x_0\in B$ and $u_0\in\pi^{-1}(x_0)$.  
Define the based $C^{\infty}$-loop spaces $\Lambda_{x_0}^{\infty}(B)$ and $\Lambda_{u_0}^{\infty}(P)$ by 
$$\begin{array}{l}
\hspace{1.5truecm}\displaystyle{\Lambda^{\infty}_{x_0}(B):=\{c\in C^{\infty}(S^1,B)\,|\,c(1)=x_0\}}\\
\displaystyle{{\rm and}\quad\,\,\,\Lambda^{\infty}_{u_0}(P):=\{\sigma\in C^{\infty}(S^1,P)\,|\,\sigma(1)=u_0\},}
\end{array}$$
respectively.  Similarly define the based $C^{\infty}$-loop group $\Lambda_e^{\infty}(G)$ at the identity element $\ve$ of $G$ 
and the based $C^{\infty}$-loop algebra $\Lambda^{\infty}_{\bf 0}(\mathfrak g)$ at the zero vector ${\bf 0}$ of $\mathfrak g$ by 
$$\begin{array}{l}
\hspace{1truecm}\displaystyle{\Lambda^{\infty}_{e}(G):=\{{\bf g}\in C^{\infty}(S^1,G)\,|\,{\bf g}(1)=e\},}\\
\displaystyle{{\rm and}\quad\,\,\Lambda^{\infty}_{{\bf 0}}(\mathfrak g):=\{u\in C^{\infty}(S^1,\mathfrak g)\,|\,u(1)={\bf 0}\},}
\end{array}$$
respectively.  
Let $\sigma$ be the horizontal lift of $c\circ z$ starting from $u_0$ with respect to $\omega_0$.  
Denote by $\pi^c:c^{\ast}P\to S^1$ the induced bundle of $P$ by $c$, which is identified with the trivial $G$-bundle 
$P_o:=S^1\times G$ over $S^1$ by $\sigma$.  
Define an immersion $\iota_c$ of the induced bundle $c^{\ast}P$ into $P$ by 
$\iota_c(z(t),u)=u\,\,\,((z(t),u)\in c^{\ast}P)$.  

\vspace{0.25truecm}

\noindent
{\bf Definition 2.1.}\ 
Define a map ${\rm hol}_c:{\mathcal A}_P^{H^s}\to G$ by 
$$(\mathcal P_{c\circ z}^{\omega}\circ(P_{c\circ z}^{\omega_0})^{-1})(u_0)=u_0\cdot{\rm hol}_c(\omega),$$
where $s$ is any non-negative integer and $\mathcal P_{c\circ z}^{\omega}$ (resp. $\mathcal P_{c\circ z}^{\omega_0}$) denotes the parallel translation along $c\circ z$ 
with respect to $\omega$ (resp. $\omega_0$).  
We call this map ${\rm hol}_c$ the {\it holonomy map along} $c$.  

\vspace{0.25truecm}

In particular, in the case of the trivial $G$-bundle $P_o:=S^1\times G$, 
the adjoint bundle ${\rm Ad}(P_o)$ is identified with the trivial $\mathfrak g$-bundle $P'_o:=S^1\times\mathfrak g$ over $S^1$, 
$\mathcal A_{P_o}^{\infty}(\approx\Omega_1^{\infty}(S^1,{\rm Ad}(P_o)))$ 
is identified with the space $C^{\infty}(S^1,\mathfrak g)$ of all $C^{\infty}$-maps of $S^1$ into $\mathfrak g$ and 
$\mathcal G_{P_o}^{\infty}$ is identified with the group $C^{\infty}(S^1,G)$ of all $C^{\infty}$-maps of $S^1$ into $G$.  
Let $s$ be a non-negative integer.  
The Hilbert space $\mathcal A_{P_o}^{H^s}(\approx\Omega_1^{H^s}(S^1,{\rm Ad}(P_o)))$ is identified with the Hilbert space $H^s([0,1],\mathfrak g)$ of 
all $H^s$-maps of $[0,1]$ into $\mathfrak g$ and $\mathcal G_{P_o}^{H^{s+1}}$ is identified with the Hilbert Lie group $H^{s+1}([0,1],G)$ of all $H^{s+1}$-maps of $[0,1]$ 
into $G$.  Here we note that, for any non-negative integer $s$, the action $\mathcal G_{P^o}^{H^{s+1}}\curvearrowright\mathcal A_{P^o}^{H^s}$ is a smooth action 
because $S^1$ is of one-dimension.  

\vspace{0.25truecm}

\noindent
{\it Remark 2.1.}\ We shall explain why $\mathcal A_{P_o}^{H^s}(\approx\Omega_1^{H^s}(S^1,{\rm Ad}(P_o)))$ is identified with the Hilbert space $H^s([0,1],\mathfrak g)$.  
The space $C^{\infty}(S^1,\mathfrak g)$ is embedded into $C^{\infty}([0,1],\mathfrak g)$ by assigning $c\circ z$ to each $c\in C^{\infty}(S^1,\mathfrak g)$.  
The orthonormal basis of the $H^s$-completion $H^s([0,1],\mathfrak g)$ of $C^{\infty}([0,1],\mathfrak g)$ is given as a family of $C^{\infty}$-loops 
in $\mathfrak g$ (having $[0,1]$ as the domain).  
Hence the Fourier expansion of each element of the $H^s$-completion $H^s([0,1],\mathfrak g)$ of $C^{\infty}([0,1],\mathfrak g)$ is given as a series consisting of 
$C^{\infty}$-loops in $\mathfrak g$.  Therefore, we see that $H^s([0,1],\mathfrak g)$ is the completion of $C^{\infty}(S^1,\mathfrak g)$.  

\vspace{0.25truecm}

As stated in Introduction, the parallel transport map for $G$ has been defined as a map of $H^0([0,1],\mathfrak g)$ onto $G$.  
In more general, we can define the parallel transport map for $G$ as the map of $H^s([0,1],\mathfrak g)$ onto $G$ similary, where $s$ is any non-negative integer.  

Let $c$ and $\sigma$ be as above.  
Then the pull-back bundle $c^{\ast}P$ of $P$ by $c$ is identified with the trivial $G$-bundle $S^1\times G$ over $S^1$ by $\sigma$ as follows.  
Define a map $\eta:S^1\times G\to c^{\ast}P$ by 
$$\eta(z(t),g):=(z(t),\sigma(t)g)\quad\,\,((t,g)\in[0,1]\times G).$$
It is clear that $\eta$ is a bundle isomorphism.  
Throughout this bundle isomorphism $\eta$, $c^{\ast}P$ is identified with $S^1\times G$.  
Similarly, a bundle isomorphism $\bar{\eta}:[0,1]\times G\to (c\circ z)^{\ast}P$ is defined by 
$$\bar{\eta}(t,g):=(t,\sigma(t)g)\quad\,\,((t,g)\in[0,1]\times G).$$
Throughout this bundle isomorphism $\bar{\eta}$, $(c\circ z)^{\ast}P$ is identified with the trivial bundle $[0,1]\times G$.  
The natural embedding $\iota_{\sigma}:[0,1]\times G\hookrightarrow P$ by 
$$\iota_{\sigma}(t,g):=\sigma(t)g\quad\,\,((t,g)\in[0,1]\times G).$$
For $\omega\in\mathcal A_P^{H^s}$, the pull-back connection $\iota_{\sigma}^{\ast}\omega$ of $(c\circ z)^{\ast}P$ is defined.  
Then we have 
$$(\iota_{\sigma}^{\ast}\omega)_{(t_0,g)}\left(\left(\frac{\partial}{\partial t}\right)_{(t_0,g)}\right)
:=\omega_{\sigma(t_0)}(\sigma'(t_0))\quad\,\,((t_0,g)\in[0,1]\times G)).$$
By the one-to-one correspondence 
$$\iota_{\sigma}^{\ast}\omega\quad\longleftrightarrow\quad t\mapsto\omega_{\sigma(t)}(\sigma'(t)) \,\,\,(t\in[0,1]),$$
$\mathcal A_{(c\circ z)^{\ast}P}^{H^s}$ is identified with the Hilbert space $H^s([0,1],\mathfrak g)$ of all $H^s$-paths in $\mathfrak g$.  

\vspace{0.25truecm}

\noindent
{\bf Definition 2.2.} Define a map $\mu_c:\mathcal A_P^{H^s}\to H^s([0,1],\mathfrak g)$ by 
$$(\mu_c(\omega))(t):=\omega_{\sigma(t)}(\sigma'(t))\quad\,\,(t\in[0,1],\,\,\,\,\omega\in\mathcal A_P^{H^s}).\leqno{(2.3)}$$
As above, $\omega_{\sigma(z(t_0))}(\sigma'(t_0))$ is equal to 
$\displaystyle{(\sigma^{\ast}\omega)_{(t_0,g)}\left(\left(\frac{\partial}{\partial t}\right)_{(t_0,g)}\right)}$.  
From this fact, we call this map $\mu_c$ the {\it pull-back connection map by $c$}.  

\vspace{0.25truecm}

\noindent
{\bf Lemma 2.1.} {\sl Among ${\rm hol}_c,\,\phi$ and $\mu_c$, the relation ${\rm hol}_c=\phi\circ\mu_c$ holds.}

\vspace{0.25truecm}

\noindent
{\it Proof.}\ Take any $\omega\in\mathcal A_P^{H^s}$.  For the simplicity, set $u:=\mu_c(\omega)$.  Let ${\bf g}_u:[0,1]\to G$ be the element of 
$H^{s+1}([0,1],G)$ satisfying ${\bf g}_u(0)=e$ and $(R_{{\bf g}_u(t)})_{\ast}^{-1}({\bf g}_u'(t))$\newline
$=u(t)$ ($t\in[0,1]$).  
Since $u(t)=\omega_{\sigma(t)}(\sigma'(t))$, ${\bf g}_u$ is the projection of the horizontal curve starting from $(0,e)$ with respect to the connection 
$\iota_{\sigma}^{\ast}\omega$ of $(c\circ z)^{\ast}P=[0,1]\times G$ onto the $G$-component.  
From this fact, we see that ${\bf g}_u(1)(=\phi(u))$ is equal to ${\rm hol}_c(\omega)$.  Hence we obtain ${\rm hol}_c(\omega)=(\phi\circ\mu_c)(\omega)$.  \qed

\vspace{0.25truecm}

By the one-to-one correspondence 
$${\bf g}\quad\longleftrightarrow\quad t\mapsto\widehat{\bf g}(t,\sigma(t))\,\,\,(t\in[0,1]),$$
$\mathcal G_{(c\circ z)^{\ast}P}^{H^{s+1}}$ is identified with the Hilbert space $H^{s+1}([0,1],G)$ of all $H^s$-paths in $G$, 
where $\widehat{\bf g}$ is the element of $\widehat{\mathcal G}_{(c\circ z)^{\ast}P}^{H^{s+1}}$ corresponding to ${\bf g}\in\mathcal G_{(c\circ z)^{\ast}P}^{H^{s+1}}$.  
Set 
$$\begin{array}{c}
\hspace{1.5truecm}\Lambda_e^{H^{s+1}}(G):=\{{\bf g}\in H^{s+1}([0,1],G)\,\vert\,{\bf g}(0)={\bf g}(1)=e\},\\
{\rm and}\quad\,\,\Lambda^{H^{s+1}}(G):=\{{\bf g}\in H^{s+1}([0,1],G)\,\vert\,{\bf g}(0)={\bf g}(1)\}.
\end{array}$$
Note that the orbit space $H^s([0,1],\mathfrak g)/\Lambda_e^{H^{s+1}}(G)$ of the subaction $\Lambda_e^{H^{s+1}}(G)\curvearrowright H^s([0,1],\mathfrak g)$ 
of the gauge action $H^{s+1}([0,1],G)\curvearrowright H^s([0,1],\mathfrak g)$ is isomorphic to 
$G$ and that the orbit map of this subaction coincides with the parallel transport map $\phi$ for $G$.  
Also, note that the orbit space $H^s([0,1],\mathfrak g)/\Lambda^{H^{s+1}}(G)$ of this subaction $\Lambda^{H^{s+1}}(G)\curvearrowright H^s([0,1],\mathfrak g)$ is isomorphic to 
$G/{\rm Ad}(G)$ and that the orbit space $H^s([0,1],\mathfrak g)/\Lambda_e^{H^{s+1}}(G)$ of this subaction $\Lambda_e^{H^{s+1}}(G)\curvearrowright H^s([0,1],\mathfrak g)$ 
is isomorphic to $G$, where ${\rm Ad}$ is the adjoint action of $G$ on oneself.  
In fact, according to the proof of Theorem 4.1 of \cite{TT}, we have 
$$\phi\circ{\bf g}=L_{{\bf g}(0)}\circ R_{{\bf g}(1)}^{-1}\circ\phi\leqno{(2.4)}$$
for any ${\bf g}\in H^{s+1}([0,1],G)$, where ${\bf g}$ in the left-hand side means the diffeomorphism of $H^s([0,1],\mathfrak g)$ onto oneself defined by the action of 
$H^{s+1}([0,1],G)$ on $H^0([0,1],\mathfrak g)$ and $L_{{\bf g}(0)}$ and $R_{{\bf g}(1)}$ are the left translation by ${\bf g}(0)$ and the right translation by ${\bf g}(1)$, 
respectively.  

\vspace{0.25truecm}

\noindent
{\bf Definition 2.3.} Define a map $\lambda_c:\mathcal G_P^{H^{s+1}}\to H^{s+1}([0,1],G)$ by 
$$\lambda_c({\bf g})(t):=\widehat{\bf g}(\sigma(t))\quad\,\,(t\in[0,1],\,\,\,\,{\bf g}\in\mathcal G_P^{H^{s+1}}).\leqno{(2.5)}$$

\vspace{0.25truecm}

The {\it based gauge transformation group} $(\mathcal G_P^{H^{s+1}})_x$ {\it at} $x\in B$ is defined by 
$$(\mathcal G_P^{H^{s+1}})_x:=\{{\bf g}\in\mathcal G_P^{H^{s+1}}\,\vert\,\widehat{\bf g}(\pi^{-1}(x)))=\{e\}\}.$$
Denote by $\mathcal M_P^{H^s}$ the muduli space $\mathcal A_P^{H^s}/\mathcal G_P^{H^{s+1}}$ and 
$\pi_{\mathcal M_P}$ the orbit map of the action $\mathcal G_P^{H^{s+1}}\curvearrowright\mathcal A_P^{H^s}$.  
Also, denote by $(\mathcal M_P^{H^s})_x$ the muduli space \newline
$\mathcal A_P^{H^s}/(\mathcal G_P^{H^{s+1}})_x$ and 
$\pi_{(\mathcal M_P)_x}$ the orbit map of the action $(\mathcal G_P^{H^{s+1}})_x\curvearrowright\mathcal A_P^{H^s}$.  
Throughout $\lambda_c$, $\mathcal G_P^{H^{s+1}}$ acts on $H^s([0,1],\mathfrak g)$.  
Also, $\mathcal G_P^{H^{s+1}}$ acts on $G$ by 
$${\bf g}\cdot g:={\rm Ad}({\bf g}(\sigma(1)))(g)\quad\,\,({\bf g}\in\mathcal G_P^{H^{s+1}},\,\,\,\,g\in G).$$
From Lemma 2.1, $(2.2),\,(2.4)$ and the defintions of $\mu_c$ and $\lambda_c$, we can derive the following facts.  

\vspace{0.25truecm}

\noindent
{\bf Lemma 2.2.} {\sl 
{\rm (i)} The pull-back connection map $\mu_c$ is $\mathcal G_P^{H^{s+1}}$-equivariant, that is, 
the following relation holds:
$$\mu_c({\bf g}\cdot\omega)=\lambda_c({\bf g})\cdot\mu_c(\omega)\quad\,\,
({\bf g}\in\mathcal G_P^{H^{s+1}},\,\,\omega\in\mathcal A_P^{H^s}).\leqno{(2.6)}$$

{\rm (ii)} The holonomy map ${\rm hol}_c$ is $\mathcal G_P^{H^{s+1}}$-equivariant, that is, 
the following relation holds:
$${\rm hol}_c({\bf g}\cdot\omega)={\rm Ad}({\bf g}(\sigma(1)))({\rm hol}_c(\omega))\quad\,\,
({\bf g}\in\mathcal G_P^{H^{s+1}},\,\,\omega\in\mathcal A_P^{H^s}).\leqno{(2.7)}$$

{\rm(iii)} Set $x_0:=c(0)$.  The pull-back connection map $\mu_c$ maps 
$(\mathcal G_P^{H^{s+1}})_{x_0}$-orbits in $\mathcal A_P^{H^s}$ to $\Lambda_e^{H^{s+1}}(G)$-orbits in $H^s([0,1],\mathfrak g)$ 
and hence there uniquely exists the map $\overline{\mu}_c$ between the orbit spaces 
$(\mathcal M_P^{H^s})_{x_0}$ and $H^s([0,1],\mathfrak g)/\Lambda_e^{H^{s+1}}(G)$\newline
$(\approx G)$ satisfying 
$$\overline{\mu}_c\circ\pi_{(\mathcal M_P)_{x_0}}=\phi\circ\mu_c(={\rm hol}_c)\leqno{(2.8)}$$
(see Diagram 2.1).  

{\rm(iv)} The pull-back connection map $\mu_c$ maps $\mathcal G_P^{H^{s+1}}$-orbits in $\mathcal A_P^{H^s}$ to \newline
$\Lambda^{H^{s+1}}(G)$-orbits in $H^s([0,1],\mathfrak g)$ 
and hence there uniquely exists the map $\overline{\overline{\mu}}_c$ between the orbit spaces $\mathcal M_P^{H^s}$ and 
$H^s([0,1],\mathfrak g)/\Lambda^{H^{s+1}}(G)(\approx G/{\rm Ad}(G))$ satisfying 
$$\overline{\overline{\mu}}_c\circ\pi_{\mathcal M_P}=\pi_{{\rm Ad}}\circ\phi\circ\mu_c(=\pi_{\rm Ad}\circ{\rm hol}_c=:\overline{\rm hol}_c)\leqno{(2.9)}$$
(see Diagram 2.2), where $\pi_{\rm Ad}$ denotes the natural projection of $G$ onto \newline
$G/{\rm Ad}(G)$.  
}

\vspace{0.25truecm}

\noindent
{\it Proof.}\ First we shall show $(2.6)$.  From $(2.2)$ and the definition of $\mu_c$, we have 
\begin{align*}
&\mu_c({\bf g}\cdot\omega)(t)=({\bf g}\cdot\omega)_{\sigma(t)}(\sigma'(t))\\
=&({\rm Ad}({\bf g}(\sigma(t))(\omega_{\sigma(t)}(\sigma'(t)))-(R_{\widehat{\bf g}(\sigma(t))})_{\ast}^{-1}((\widehat{\bf g}\circ\sigma)'(t)).
\end{align*}
On the other hand, from $(2.2)$ and the definitions of $\lambda_c$ and $\mu_c$, we have 
\begin{align*}
&(\lambda_c({\bf g})\cdot\mu_c(\omega))(t)=\left(\,\,(\widehat{\bf g}\circ\sigma)\cdot(\omega\circ\sigma')\,\,\right)(t)\\
=&{\rm Ad}({\bf g}(\sigma(t)))(\omega_{\sigma(t)}(\sigma'(t)))-(R_{{\bf g}(\sigma(t)})_{\ast}^{-1}((\widehat{\bf g}\circ\sigma)'(t)).
\end{align*}
Therefore we obtain $\mu_c({\bf g}\cdot\omega)(t)=(\lambda_c({\bf g})\cdot\mu_c(\omega))(t)$.  

Next we shall show $(2.7)$.  From ${\rm hol}_c=\phi\circ\mu_c$ (by Lemma 2.1) and $(2.4)$, we have 
\begin{align*}
&{\rm hol}_c({\bf g}\cdot\omega)=\phi(\,({\bf g}\cdot\omega)\circ\sigma'\,)\\
=&\phi\left(\,\,({\rm Ad}({\bf g}\circ\sigma)(\omega\circ\sigma')-(R_{\widehat{\bf g}\circ\sigma})_{\ast}^{-1}\circ(\widehat{\bf g}\circ\sigma)'
\,\,\right)\\
=&\phi(\,({\bf g}\circ\sigma)\cdot(\omega\circ\sigma')\\
=&\left(L_{{\bf g}(\sigma(z(0))}\circ R_{{\bf g}(\sigma(z(1))}^{-1}\right)(\phi(\omega\circ\sigma'))\\
=&{\rm Ad}({\bf g}(\sigma(1)))\,({\rm hol}_c(\omega)).
\end{align*}

Next we shall show the statements (iii) and (iv).  It is clear that \newline
$\lambda_c((\mathcal G_P^{H^{s+1}})_{x_0})=\Lambda_e^{H^{s+1}}(G)$ and 
$\lambda_c(\mathcal G_P^{H^{s+1}})=\Lambda^{H^{s+1}}(G)$ hold.  
From these facts and $(2.6)$, we can derive the statements (iii) and (iv) directly.  \qed

\vspace{0.5truecm}

{\small 
\centerline{
\unitlength 0.1in
\begin{picture}( 52.7000, 15.0500)(-16.4000,-21.2500)
\put(18.9000,-13.5000){\makebox(0,0)[rt]{$\mathcal A_P^{H^s}$}}%
\put(31.7000,-13.5000){\makebox(0,0)[lt]{$H^s([0,1],\mathfrak g)$}}%
\put(12.4000,-9.1000){\makebox(0,0)[rb]{$(\mathcal G_P^{H^{s+1}})_{x_0}$}}%
\put(30.4000,-9.1000){\makebox(0,0)[rb]{$\Lambda_e^{H^{s+1}}(G)$}}%
%
\special{pn 8}%
\special{pa 1440 850}%
\special{pa 2330 850}%
\special{fp}%
\special{sh 1}%
\special{pa 2330 850}%
\special{pa 2264 830}%
\special{pa 2278 850}%
\special{pa 2264 870}%
\special{pa 2330 850}%
\special{fp}%
\put(17.8000,-7.9000){\makebox(0,0)[lb]{$\lambda_c$}}%
%
\special{pn 8}%
\special{pa 2050 1430}%
\special{pa 3000 1430}%
\special{fp}%
\special{sh 1}%
\special{pa 3000 1430}%
\special{pa 2934 1410}%
\special{pa 2948 1430}%
\special{pa 2934 1450}%
\special{pa 3000 1430}%
\special{fp}%
\put(23.7000,-13.8000){\makebox(0,0)[lb]{$\mu_c$}}%
%
\special{pn 13}%
\special{ar 1100 1310 590 450  5.0767644 5.9136541}%
%
\special{pn 13}%
\special{pa 1660 1150}%
\special{pa 1700 1230}%
\special{fp}%
\special{sh 1}%
\special{pa 1700 1230}%
\special{pa 1688 1162}%
\special{pa 1676 1182}%
\special{pa 1652 1180}%
\special{pa 1700 1230}%
\special{fp}%
%
\special{pn 13}%
\special{ar 2900 1320 590 450  5.0767644 5.9136541}%
%
\special{pn 13}%
\special{pa 3460 1160}%
\special{pa 3500 1240}%
\special{fp}%
\special{sh 1}%
\special{pa 3500 1240}%
\special{pa 3488 1172}%
\special{pa 3476 1192}%
\special{pa 3452 1190}%
\special{pa 3500 1240}%
\special{fp}%
%
\special{pn 8}%
\special{pa 1710 1610}%
\special{pa 1710 1980}%
\special{fp}%
\special{sh 1}%
\special{pa 1710 1980}%
\special{pa 1730 1914}%
\special{pa 1710 1928}%
\special{pa 1690 1914}%
\special{pa 1710 1980}%
\special{fp}%
%
\special{pn 8}%
\special{pa 3550 1580}%
\special{pa 3550 1950}%
\special{fp}%
\special{sh 1}%
\special{pa 3550 1950}%
\special{pa 3570 1884}%
\special{pa 3550 1898}%
\special{pa 3530 1884}%
\special{pa 3550 1950}%
\special{fp}%
\put(19.6000,-20.3000){\makebox(0,0)[rt]{$(\mathcal M_P^{H^s})_{x_0}$}}%
\put(34.9000,-20.4000){\makebox(0,0)[lt]{$G$}}%
%
\special{pn 8}%
\special{pa 2150 2120}%
\special{pa 3300 2120}%
\special{fp}%
\special{sh 1}%
\special{pa 3300 2120}%
\special{pa 3234 2100}%
\special{pa 3248 2120}%
\special{pa 3234 2140}%
\special{pa 3300 2120}%
\special{fp}%
\put(25.8000,-20.9000){\makebox(0,0)[lb]{$\overline{\mu}_c$}}%
%
\special{pn 8}%
\special{pa 2020 1560}%
\special{pa 3270 1970}%
\special{fp}%
\special{sh 1}%
\special{pa 3270 1970}%
\special{pa 3214 1930}%
\special{pa 3220 1954}%
\special{pa 3200 1968}%
\special{pa 3270 1970}%
\special{fp}%
\put(26.1000,-17.3000){\makebox(0,0)[lb]{{\small ${\rm hol}_c$}}}%
\put(36.3000,-17.0000){\makebox(0,0)[lt]{$\phi$}}%
%
\special{pn 8}%
\special{ar 2132 1832 72 72  0.6528466 5.7088805}%
%
\special{pn 8}%
\special{pa 2192 1792}%
\special{pa 2212 1842}%
\special{fp}%
\special{sh 1}%
\special{pa 2212 1842}%
\special{pa 2206 1774}%
\special{pa 2192 1792}%
\special{pa 2170 1788}%
\special{pa 2212 1842}%
\special{fp}%
\put(16.1000,-18.4000){\makebox(0,0)[rb]{$\pi_{(\mathcal M_P)_{x_0}}$}}%
%
\special{pn 8}%
\special{ar 3232 1692 72 72  0.6528466 5.7088805}%
%
\special{pn 8}%
\special{pa 3292 1652}%
\special{pa 3312 1702}%
\special{fp}%
\special{sh 1}%
\special{pa 3312 1702}%
\special{pa 3306 1634}%
\special{pa 3292 1652}%
\special{pa 3270 1648}%
\special{pa 3312 1702}%
\special{fp}%
\end{picture}%
\hspace{8truecm}}
}

\vspace{0.35truecm}

\centerline{{\bf Diagram 2.1$\,\,:\,\,$ The map $\overline{\mu}_c$ induced from $\mu_c$}}

\vspace{0.5truecm}

{\small 
\centerline{
\unitlength 0.1in
\begin{picture}( 45.5000, 15.0500)( -9.2000,-21.2500)
\put(18.9000,-13.5000){\makebox(0,0)[rt]{$\mathcal A_P^{H^s}$}}%
\put(31.7000,-13.5000){\makebox(0,0)[lt]{$H^s([0,1],\mathfrak g)$}}%
\put(12.4000,-9.1000){\makebox(0,0)[rb]{$\mathcal G_P^{H^{s+1}}$}}%
\put(30.4000,-9.1000){\makebox(0,0)[rb]{$\Lambda^{H^{s+1}}(G)$}}%
%
\special{pn 8}%
\special{pa 1440 850}%
\special{pa 2330 850}%
\special{fp}%
\special{sh 1}%
\special{pa 2330 850}%
\special{pa 2264 830}%
\special{pa 2278 850}%
\special{pa 2264 870}%
\special{pa 2330 850}%
\special{fp}%
\put(17.8000,-7.9000){\makebox(0,0)[lb]{$\lambda_c$}}%
%
\special{pn 8}%
\special{pa 2050 1430}%
\special{pa 3000 1430}%
\special{fp}%
\special{sh 1}%
\special{pa 3000 1430}%
\special{pa 2934 1410}%
\special{pa 2948 1430}%
\special{pa 2934 1450}%
\special{pa 3000 1430}%
\special{fp}%
\put(23.7000,-13.8000){\makebox(0,0)[lb]{$\mu_c$}}%
%
\special{pn 13}%
\special{ar 1100 1310 590 450  5.0767644 5.9136541}%
%
\special{pn 13}%
\special{pa 1660 1150}%
\special{pa 1700 1230}%
\special{fp}%
\special{sh 1}%
\special{pa 1700 1230}%
\special{pa 1688 1162}%
\special{pa 1676 1182}%
\special{pa 1652 1180}%
\special{pa 1700 1230}%
\special{fp}%
%
\special{pn 13}%
\special{ar 2900 1320 590 450  5.0767644 5.9136541}%
%
\special{pn 13}%
\special{pa 3460 1160}%
\special{pa 3500 1240}%
\special{fp}%
\special{sh 1}%
\special{pa 3500 1240}%
\special{pa 3488 1172}%
\special{pa 3476 1192}%
\special{pa 3452 1190}%
\special{pa 3500 1240}%
\special{fp}%
%
\special{pn 8}%
\special{pa 1710 1610}%
\special{pa 1710 1980}%
\special{fp}%
\special{sh 1}%
\special{pa 1710 1980}%
\special{pa 1730 1914}%
\special{pa 1710 1928}%
\special{pa 1690 1914}%
\special{pa 1710 1980}%
\special{fp}%
%
\special{pn 8}%
\special{pa 3550 1580}%
\special{pa 3550 1950}%
\special{fp}%
\special{sh 1}%
\special{pa 3550 1950}%
\special{pa 3570 1884}%
\special{pa 3550 1898}%
\special{pa 3530 1884}%
\special{pa 3550 1950}%
\special{fp}%
\put(19.6000,-20.3000){\makebox(0,0)[rt]{$\mathcal M_P^{H^s}$}}%
\put(34.9000,-20.4000){\makebox(0,0)[lt]{$G/{\rm Ad}(G)$}}%
%
\special{pn 8}%
\special{pa 2150 2120}%
\special{pa 3300 2120}%
\special{fp}%
\special{sh 1}%
\special{pa 3300 2120}%
\special{pa 3234 2100}%
\special{pa 3248 2120}%
\special{pa 3234 2140}%
\special{pa 3300 2120}%
\special{fp}%
\put(25.8000,-20.9000){\makebox(0,0)[lb]{$\overline{\overline{\mu}}_c$}}%
%
\special{pn 8}%
\special{pa 2020 1560}%
\special{pa 3270 1970}%
\special{fp}%
\special{sh 1}%
\special{pa 3270 1970}%
\special{pa 3214 1930}%
\special{pa 3220 1954}%
\special{pa 3200 1968}%
\special{pa 3270 1970}%
\special{fp}%
\put(26.1000,-17.3000){\makebox(0,0)[lb]{{\small $\overline{{\rm hol}}_c$}}}%
\put(36.3000,-17.0000){\makebox(0,0)[lt]{$\pi_{\rm Ad}\circ\phi$}}%
%
\special{pn 8}%
\special{ar 2132 1832 72 72  0.6528466 5.7088805}%
%
\special{pn 8}%
\special{pa 2192 1792}%
\special{pa 2212 1842}%
\special{fp}%
\special{sh 1}%
\special{pa 2212 1842}%
\special{pa 2206 1774}%
\special{pa 2192 1792}%
\special{pa 2170 1788}%
\special{pa 2212 1842}%
\special{fp}%
\put(16.1000,-18.4000){\makebox(0,0)[rb]{$\pi_{\mathcal M_P}$}}%
%
\special{pn 8}%
\special{ar 3232 1692 72 72  0.6528466 5.7088805}%
%
\special{pn 8}%
\special{pa 3292 1652}%
\special{pa 3312 1702}%
\special{fp}%
\special{sh 1}%
\special{pa 3312 1702}%
\special{pa 3306 1634}%
\special{pa 3292 1652}%
\special{pa 3270 1648}%
\special{pa 3312 1702}%
\special{fp}%
\end{picture}%
\hspace{6truecm}}
}

\vspace{0.35truecm}

\centerline{{\bf Diagram 2.2$\,\,:\,\,$ The map $\overline{\overline{\mu}}_c$ induced from $\mu_c$}}

\section{Proof of Theorem A} 
In this section, we shall prove Theorem A stated in Introduction.  
Fix $\omega_0\in\mathcal A_P^{\infty}$.  
Let $\mathcal V^{\mathcal A}$ be the vertical distribution on $\mathcal A_P^{H^s}$ with respect to $\pi_{\mathcal M_P}$, that is, 
$$\mathcal V^{\mathcal A}_{\omega}:={\rm Ker}\,(\pi_{\mathcal M_P})_{\ast\omega}\quad\,\,(\omega\in \mathcal A_P^{H^s})$$
and $\mathcal H^{\mathcal A}$ the horizontal distribution on $\mathcal A_P^{H^s}$ with respect to $\langle\,\,,\,\,\rangle_s^{\omega_0}$, that is, 
$$\mathcal H^{\mathcal A}_{\omega}:=\{\vv\in T_{\omega}\mathcal A_P^{H^s}\,\vert\,\langle\vv,\vw\rangle_s^{\omega_0}=0\,\,\,\,
(\forall\,\vw\in\mathcal V^{\mathcal A}_{\omega})\}.$$
Also, let $\mathcal V^P$ be the vertical distribution on $P$, that is, 
$$\mathcal V^P_u:={\rm Ker}\,\pi_{\ast u}\quad\,\,(u\in P)$$
and $\mathcal H^P$ the horizontal distribution on $P$ with respect to $\omega_0$, that is, $\mathcal H^P_u:={\rm Ker}\,(\omega_0)_u$ ($u\in P$).  
For $u\in P,\,\vv\in T_uP\setminus\mathcal V_u$ and $\xi\in\mathfrak g$, define a linear map $\eta_{\vv,\xi}:T_uP\to\mathfrak g$ by 
$$\eta_{\vv,\xi}(\vw):=\left\{
\begin{array}{ll}
\xi & (\vw=\vv)\\
{\bf 0} & (\vw\in\mathcal V^P_u\oplus(\mathcal H^P_u\cap{\rm Span}\{\vv\}^{\perp})).
\end{array}\right.$$
Define $\delta_{\eta_{\vv,\xi}}\in\Omega_{\mathcal T,1}^{H^s}(P,\mathfrak g)$ by 
$$\langle\delta_{\eta_{\vv,\xi}},A\rangle_0=\langle\widehat{\eta}_{\vv,\xi},\widehat A_{\pi(u)}\rangle_{B,\mathfrak g}
\quad(\forall\,A\in\Omega_{\mathcal T,1}^{H^0}(P,\mathfrak g)),$$
where $\widehat{\eta}_{\vv,\xi}$ is an element of $T_{\pi(u)}^{\ast}B\otimes{\rm Ad}(P)_{\pi(u)}$ 
corresponding to $\eta_{\vv,\xi}$.  

\vspace{0.25truecm}

\noindent
{\bf Lemma 3.1.}\ {\sl Let $c$ and $\sigma$ be as in the statement of Theorem A.  
Then we have 
$$\left({\rm Ker}(\mu_c)_{\ast\omega}\right)^{\perp}
=\left\{\left.\int_0^1\delta_{\eta_{\sigma'(t),\xi(t)}}\,dt\,\,\right\vert\,\,\xi\in H^s([0,1],\mathfrak g)\right\}.$$
}

\vspace{0.25truecm}

\noindent
{\it Proof.}\ \ First we note that $\mu_c$ is linear.  Since 
$${\rm Ker}(\mu_c)_{\ast\omega}=\left\{\left.A\in\Omega_{\mathcal T,1}^{H^s}(P,\mathfrak g)\,\,\right\vert\,\,A_{\sigma(t)}(\sigma'(t))
={\bf 0}\quad(\forall\,t\in[0,1])\right\},$$
$A\in\Omega_{\mathcal T,1}^{H^s}(P,\mathfrak g)$ belongs to ${\rm Ker}\,((\mu_c)_{\ast\omega})^{\perp}$ if and only if $A$ vanishes over 
$$\left(TP\setminus(TP\vert_{\sigma(S^1)})\right)\,\cup\,\left(\mathop{\amalg}_{t\in[0,1]}
(\mathcal H_{\sigma(t)}^P\cap{\rm Span}\{\sigma'(t)\}^{\perp})\right),$$
that is, $A$ is expressed as $\displaystyle{A=\int_0^1\delta_{\eta_{\sigma'(t),\xi(t)}}\,dt}$ for some $\xi\in H^s([0,1],\mathfrak g)$.  
Therefore, we obtain the desired relation.  \qed

\vspace{0.25truecm}

Set 
$$\begin{array}{c}
\hspace{0.4truecm}\Lambda^{H^s}(\mathfrak g):=\{u\in H^s([0,1],\mathfrak g)\,\vert\,u(0)=u(1)\}\\
{\rm and}\quad\,\,\Lambda_{\bf 0}^{H^s}(\mathfrak g):=\{u\in H^s([0,1],\mathfrak g)\,\vert\,u(0)=u(1)={\bf 0}\}.
\end{array}$$
Let $g_s^o$ be the Riemannian metric of $H^s([0,1],\mathfrak g)(\approx\mathcal A_{P^o}^{H^s})$ defined in similar to $g_s$.  

\vspace{0.25truecm}

\noindent
{\bf Proposition 3.2.}\ {\sl For any non-negative integer $s$, the pull-back connection map 
$$\mu_c:(\mathcal A_P^{H^s},g_s)\to(H^s([0,1],\mathfrak g),g_s^o)$$
along $c$ is a homothetic submersion of coefficient $a$ onto the dense linear subspace $\Lambda^{H^s}(\mathfrak g)$ of $H^s([0,1],\mathfrak g)$.  
Also, if $s=0$, then the fibres of $\mu_c$ are totally geodesic.}

\vspace{0.25truecm}

\noindent
{\it Proof.}\ \ 
In this proof, we abbreviate $\widehat A\in\Gamma^{H^s}(T^{\ast}B\otimes{\rm Ad}(P))$ corresponding to $A\in\Omega_{\mathcal T,1}^{H^s}(P,\mathfrak g)$ as $A$ and 
$\widehat A_{\pi(\bullet)}\in T^{\ast}_{\cdot}B\otimes{\rm Ad}(P)_{\pi(\bullet)}$ corresponding to $A_{\bullet}\in T_{\bullet}P\otimes\mathfrak g$ as $A_{\bullet}$ 
for the simplicity.  
Take $A\in\left({\rm Ker}(\mu_c)_{\ast\omega}\right)^{\perp}$.  According to Lemma 3.1, $A$ is expressed as 
$\displaystyle{A:=\int_0^1\delta_{\eta_{\sigma'(t),\xi(t)}}\,dt}$ for some $\xi\in H^s([0,1],\mathfrak g)$.  
Then we have 
$$\begin{array}{l}
\hspace{0.5truecm}\displaystyle{(g_s)_{\omega}(A,A)=\langle A,\square_{\omega}^s(A)\rangle_0}\\
\displaystyle{=\langle\int_0^1\delta_{\eta_{\sigma'(\bar t),\xi(\bar t)}}\,d\bar t,\,\,
\square_{\omega}^s\left(\int_0^1\delta_{\eta_{\sigma'(\hat t),\xi(\hat t)}}\,d\hat t\right)\rangle_0}\\
\displaystyle{=\int_0^1\int_0^1\langle\eta_{\sigma'(\bar t),\xi(\bar t)},\,\,
\square_{\omega}^s\left(\delta_{\eta_{\sigma'(\hat t),\xi(\hat t)}}\right)_{c(z(\bar t))}\rangle_{B,\mathfrak g}\,d\bar td\widehat t}\\
\displaystyle{=\int_0^1\langle\eta_{\sigma'(t),\xi(t)},\,\,
\square_{\omega}^s\left(\delta_{\eta_{\sigma'(t),\xi(t)}}\right)_{c(z(t))}\rangle_{B,\mathfrak g}\,dt}\\
\displaystyle{=\frac{1}{a^2}\,\int_0^1\langle\eta_{\sigma'(t),\xi(t)}((c\circ z)'(t)),\,\,}\\
\hspace{1.5truecm}\displaystyle{\square_{\omega}^s\left(\delta_{\eta_{\sigma'(t),\xi(t)}}\right)_{c(z(t))}
((c\circ z)'(t))\rangle_{{\rm Ad}(P)_{c(z(t))}}\,dt}\\
\displaystyle{=\frac{1}{a^2}\,\int_0^1\langle\xi(t),\,\,
\left(\square_{\omega}^s\left(\delta_{\eta_{\sigma'(t),\xi(t)}}\right)\right)_{\sigma(t)}(\sigma'(t))
\rangle_{\mathfrak g}\,dt}
\end{array}\leqno{(3.1)}$$
\newpage
$$\begin{array}{l}
\displaystyle{=\frac{1}{a^2}\,\int_0^1\langle\xi(t),\,\,
\sigma^{\ast}\left(\square_{\omega}^s\left(\delta_{\eta_{\sigma'(t),\xi(t)}}\right)\right)_t\left(\frac{d}{dt}\right)\rangle_{\mathfrak g}\,dt}\\
\displaystyle{=\frac{1}{a^2}\,\langle\xi,\,\square_{\mu_c(\omega)}^s(\xi)\rangle_0,}
\end{array}
$$
where 
$\langle\,\,,\,\,\rangle_{{\rm Ad}(P)_{c(z(t))}}$ denotes the fibre metric of the bundle ${\rm Ad}(P)$ at $c(z(t))$.  
On the other hand, we have 
\begin{align*}
((\mu_c)_{\ast\omega}(A))(t)&=\int_0^1(\delta_{\eta_{\sigma'(\bar t),\xi(\bar t)}})_{\sigma(t)}(\sigma'(t))\,d\bar t\\
&=\eta_{\sigma'(t),\xi(t)}(\sigma'(t))=\xi(t)
\end{align*}
and hence 
$$\begin{array}{l}
\hspace{0.5truecm}\displaystyle{(g_s^o)_{\mu_c(\omega)}((\mu_c)_{\ast\omega}(A),(\mu_c)_{\ast\omega}(A))}\\
\displaystyle{=\langle(\mu_c)_{\ast\omega}(A),(\mu_c)_{\ast\omega}(A)\rangle_s^{\mu_c(\omega)}}\\
\displaystyle{=\langle\xi,\xi\rangle_s^{\mu_c(\omega)}=\langle\xi,\square_{\mu_c(\omega)}^s(\xi)\rangle_0.}
\end{array}\leqno{(3.2)}$$
From $(3.1)$ and $(3.2)$, we obtain 
$$(g_s^o)_{\mu_c(\omega)}((\mu_c)_{\ast\omega}(A),(\mu_c)_{\ast\omega}(A))=a^2\,(g_s)_{\omega}(A,A).$$
Also, it is clear that $\mu_c(\mathcal A_P^{H^s})=\Lambda^{H^s}(\mathfrak g)$ holds.  
These facts imply that $\mu_c$ is a homothetic submersion of coefficient $a$ of $(\mathcal A_P^{H^s},g_s)$ 
onto a dense linear subspace $\Lambda^{H^s}(\mathfrak g)$ of $(H^s([0,1],\mathfrak g),g_s^o)$.  
In particular, in the case of $s=0$, $(\mathcal A_P^{H^0},g_0)$ and $(H^0([0,1],\mathfrak g),g_0^o)$ are Hilbert spaces 
and $\mu_c$ is a linear map between these Hilbert spaces.  
Therefore the fibres of $\mu_c$ are affine subspaces of $\mathcal A_P^{H^0}$ and hence they are totally geodesic in $(\mathcal A_P^{H^0},g_0)$.  
This completes the proof.  \qed

\vspace{0.25truecm}

For the parallel transport map, we have the following fact.  

\vspace{0.25truecm}

\noindent
{\bf Proposition 3.3} {\sl For any non-negative integer $s$, the parallel transport map $\phi:(H^s([0,1],\mathfrak g),g_s^o)\to (G,g_G)$ is a Riemannian submersion.}

\vspace{0.25truecm}

\noindent
{\it Proof.}\ The gauge action $H^{s+1}([0,1],G)\curvearrowright(H^s([0,1],\mathfrak g),g_s^o)$ is isometric (and free) and 
$\phi$ is the orbit map of this subaction $\Lambda_e^{H^{s+1}}(G)\curvearrowright H^s([0,1],\mathfrak g)$, 
there exists the metric $\bar g$ of $G$ such that $\phi:(H^s([0,1],\mathfrak g),g_s^o)\to(G,\bar g)$ is a Riemannian submersion.  
It is easy to show that $\bar g$ is equal to the bi-invariant metric $g_G$.  Therefore we obtain the statement of this proposition.  \qed

\vspace{0.25truecm}

Let $K$ be a symmetric subgroup of $G$, that is, a closed subgroup of $G$ such that $({\rm Fix}\,\theta)_0\subset K\subset{\rm Fix}\,\theta$ holds for some involution of $G$, 
where ${\rm Fix}\,\theta$ denotes the fixed point set of $\theta$ and $({\rm Fix}\,\theta)_0$ denotes the identity component of ${\rm Fix}\,\theta$.  
Since the natural $K$-action on $(G,g_G)$ is isometric, there exists a unique Riemannian metric $g_{G/K}$ on $G/K$ such that 
the natural projection $\pi_{G/K}:(G,g_G)\to(G/K,g_{G/K})$ is a Riemannian submersion.  The Riemannian manifold $(G/K,g_{G/K})$ is a symmetric space of compact type.  
Let $\mathfrak g=\mathfrak k\oplus\mathfrak p$ be the canonical decomposition associated to the symmetric pair $(G,K)$.  
The space $\mathfrak p$ is identified with the tangent space $T_{eK}(G/K)$ through $(\pi_{G/K})_{\ast e}\vert_{\mathfrak p}$.  
Take a maximal abelian subspace $\mathfrak a$ of $\mathfrak p$ and let 
$$\mathfrak p=\mathfrak a\oplus\left(\mathop{\oplus}_{\alpha\in\triangle_+}\mathfrak p_{\alpha}\right)\quad\,\,{\rm and}\quad\,\,
\mathfrak k=\mathfrak z_{\mathfrak k}(\mathfrak a)\oplus\left(\mathop{\oplus}_{\alpha\in\triangle_+}\mathfrak k_{\alpha}\right)$$
be the rooot space decomposition of $\mathfrak p$ (resp. $\mathfrak k$) with respect to $\mathfrak a$ ($\triangle_+\,:\,$ the positive root system), where 
$\mathfrak z_{\mathfrak k}(\mathfrak a)$ is the centralizer of $\mathfrak a$ in $\mathfrak k$ and 
$\mathfrak p_{\alpha}$ (resp. $\mathfrak k_{\alpha}$) is the root space for $\alpha$, that is, 
$$\begin{array}{l}
\displaystyle{\mathfrak p_{\alpha}=\{\vw\in\mathfrak p\,\vert\,{\rm ad}(\vv)^2(\vw)=-\alpha(\vv)^2\,\vw\,\,\,\,(\forall\,\vv\in\mathfrak a)\},}\\
\displaystyle{\mathfrak k_{\alpha}=\{\vw\in\mathfrak k\,\vert\,{\rm ad}(\vv)^2(\vw)=-\alpha(\vv)^2\,\vw\,\,\,\,(\forall\,\vv\in\mathfrak a)\}.}
\end{array}$$
Set $\widetilde{\mathfrak a}:=\mathfrak a\oplus\mathfrak z_{\mathfrak k}(\mathfrak a)$.  
For $\vv\in\mathfrak g$, denote by $\widehat{\vv}$ the element of $H^0([0,1],\mathfrak g)$ defined as the constant map at $\vv$.  
Let $\ve^{\alpha}$ be a unit vector belonging to $\mathfrak p_{\alpha}$ and $\ve^0$ a unit vector belonging to $\widetilde{\mathfrak a}$.  
Define $l^i_{\ve^{\alpha},j}\in H^0([0,1],\mathfrak g)$ ($i=1,2,\,\,\,j\in\mathbb Z\setminus\{0\}$) by 
$$\begin{array}{l}
\displaystyle{l^1_{\ve^{\alpha},j}(z(t)):=
\ve^{\alpha}\,\cos(2j\pi t)-(\ve^{\alpha}_{\mathfrak k})\sin(2j\pi t),}\\
\displaystyle{l^2_{\ve^{\alpha},j}(z(t)):=
\ve^{\alpha}\,\sin(2j\pi t)+(\ve^{\alpha}_{\mathfrak k})\cos(2j\pi t)}
\end{array}$$
($t\in[0,1]$), where $\ve^{\alpha}_{\mathfrak k}$ is the element of $\mathfrak k_{\alpha}$ such that 
$${\rm ad}(\vv)(\ve^{\alpha})=\alpha(\vv)\,\ve^{\alpha}_{\mathfrak k}\,\,\,{\rm and}\,\,\,{\rm ad}(\vv)(\ve^{\alpha}_{\mathfrak k})=-\alpha(\vv)\,\ve^{\alpha}$$
for any $\vv\in\mathfrak a$.  
Also, define $l^i_{\ve^0,j}\in H^0([0,1],\mathfrak g)$ ($i=1,2,\,\,\,j\in\mathbb N\setminus\{0\}$) by 
$$\begin{array}{l}
\displaystyle{l^1_{\ve^0,j}(z(t)):=
\cos(2j\pi t)\cdot\ve^0\,\,\,\,(t\in[0,1]),}\\
\displaystyle{l^2_{\ve^0,j}(z(t)):=
\sin(2j\pi t)\cdot\ve^0\,\,\,\,(t\in[0,1]).}
\end{array}$$
Then, from the definition of the $L^2$-inner product $\langle\,\,,\,\,\rangle_0(=g_0^o)$ of $H^0([0,1],\mathfrak g)$, we can show the following facts directly.  

\vspace{0.25truecm}

\noindent
{\bf Lemma 3.4.}\ {\sl Let $\{\ve^{\alpha}_j\}_{j=1}^{m_{\alpha}}$ be an orthonormal basis of $\mathfrak p_{\alpha}$ and $\{\ve^0_j\}_{j=1}^{m_0}$ that of 
$\widetilde{\mathfrak a}$.  Then the statements {\rm (i)}, {\rm (ii)}  and {\rm (iii)} hold:

{\rm (i)}\ \ The system 
$$\begin{array}{c}
\displaystyle{\{\widehat{\ve^0_j}\,\vert\,1\leq j\leq m_0\}\amalg
\left(\mathop{\amalg}_{\alpha\in\triangle_+}\{\widehat{\ve^{\alpha}_j}\,\vert\,1\leq j\leq m_{\alpha}\}\right)}\\
\displaystyle{\amalg\{l^i_{\ve^0_j,k}\,\vert\,i=1,2,\,\,\,1\leq j\leq m_0,\,\,k\in\mathbb N\setminus\{0\}\}}\\
\displaystyle{\amalg\left(\mathop{\amalg}_{\alpha\in\triangle_+}\{l^i_{\ve^{\alpha}_j,k}\,\vert\,i=1,2,\,\,\,1\leq j\leq m_{\alpha},\,\,k\in\mathbb Z\setminus\{0\}\}\right)}
\end{array}$$
is an orthonormal basis of $(H^0([0,1],\mathfrak g),g_0^o)$.  

{\rm (ii)}\ \ The system 
$$\{\widehat{\ve^0_j}\,\vert\,1\leq j\leq m_0\}\amalg
\left(\mathop{\amalg}_{\alpha\in\triangle_+}\{\widehat{\ve^{\alpha}_j}\,\vert\,1\leq j\leq m_{\alpha}\}\right)$$
is an orthonormal basis of the horizontal space (which is denoted by $\mathcal H^{\phi}_{\hat 0}$) of $\phi$ at $\hat 0$.  

{\rm (iii)}\ \ The system 
$$\begin{array}{c}
\displaystyle{\{l^i_{\ve^0_j,k}\,\vert\,i=1,2,\,\,\,1\leq j\leq m_0,\,\,k\in\mathbb N\setminus\{0\}\}}\\
\displaystyle{\amalg\left(\mathop{\amalg}_{\alpha\in\triangle_+}\{l^i_{\ve^{\alpha}_j,k}\,\vert\,i=1,2,\,\,\,1\leq j\leq m_{\alpha},\,\,k\in\mathbb Z\setminus\{0\}\}\right)}
\end{array}$$
is an orthonormal basis of the vertical space (which is denoted by $\mathcal V^{\phi}_{\hat 0}$) of $\phi$ at $\hat 0$.  
}

\vspace{0.25truecm}

By using $(2.4)$, we can show that the action $H^1([0,1],G)\curvearrowright(H^0([0,1],\mathfrak g),$\newline
$g_0^o)$ is isometric and transitive 
and each element of $H^1([0,1],G)$ maps the fibres of $\phi$ to them.  Hence the fibres of $\phi$ are congruent to one another.  
We suffice to focus attention on the fibre $\phi^{-1}(e)$ (which passes through $\hat 0$).  
Denote by $A^F$ the shape tensor of $M:=\phi^{-1}(e)$.  

\vspace{0.25truecm}

\noindent
{\bf Lemma 3.5.}\ {\sl Let $\vv\in T_eG(=\mathfrak g)$ and $\widehat{\vv}(\in H^0([0,1],\mathfrak g))$ the constant path at $\vv$.  Then the following relations hold:
$$\begin{array}{c}
\displaystyle{A^F_{\hat{\vv}}(l^i_{\ve^0_j,k})={\bf 0}\quad(i=1,2,\,\,\,j=1,\cdots,m_0,\,\,\,k\in\mathbb N\setminus\{0\}),}\\
\displaystyle{A^F_{\hat{\vv}}(l^i_{\ve^{\alpha}_j,k})=-\frac{\alpha(\vv)}{2k\pi}\,l^i_{\ve^{\alpha}_j,k}\quad(i=1,2,\,\,j=1,\cdots,m_{\alpha},\,\,k\in\mathbb Z\setminus\{0\}).}
\end{array}$$
}

\vspace{0.25truecm}

\noindent
{\it Proof.}\ \ 
Denote by $\widetilde{\nabla}$ the Riemannian connection of the Hilbert space \newline
$(H^0([0,1],\mathfrak g),g_0^o)$.  
From the relations in the proof of Propositions 3.1 and 3.2 of \cite{K}, 
we obtain 
$$\widetilde{\nabla}_{l^i_{\ve^0_j,k}}\vv^L={\bf 0}\quad\,\,{\rm and}\quad\,\,\widetilde{\nabla}_{l^i_{\ve^{\alpha}_j,k}}\vv^L=\frac{\alpha(\vv)}{2k\pi}\,l^i_{\ve^{\alpha}_j,k}.
\leqno{(3.5)}$$
where $\vv^L$ is the horizontal lift of $\vv$.  Furthermore, from these relations, we obtain the desired relations.  \qed

\vspace{0.25truecm}

From this lemma, we can derive the following fact.  

\vspace{0.25truecm}

\noindent
{\bf Proposition 3.6.}\ {\sl The parallel transport map $\phi:(H^0([0,1],\mathfrak g),g_0^o)\to(G,g_G)$ 
is a Riemannian submersion with minimal regularizable fibres.}

\vspace{0.25truecm}

\noindent
{\it Proof.}\ From the relations in Lemma 3.5, it follows that $M=\phi^{-1}(e)$ is a minimal regularizable submanifold in 
$(H^0([0,1],\mathfrak g),g_0^o)$.  Hence so are all the fibres.  \qed

\vspace{0.25truecm}

According to Theorem 6.5 of \cite{HLO}, the following fact holds.  

\vspace{0.25truecm}

\noindent
{\bf Theorem 3.7(\cite{HLO}).}\ {\sl Let $\psi$ be a Riemannian submersion with minimal regularizable fibres of a Hilbert space $V$ onto a finite dimensional 
Riemannian manifold $N$.  Then, for a compact submanifold $M$ in $N$, the following statements are equivalent:

{\rm (i)}\ \ $M$ is equifocal;

{\rm (ii)}\ \ $\psi^{-1}(M)$ is isoparametric.
}

\vspace{0.25truecm}

By using Propositions 3.2, 3.3 and 3.6 and Theorem 3.7, we prove Theorem A.  

\vspace{0.25truecm}

\noindent
{\bf Proof of Theorem A.}\ 
Since ${\rm hol}_c=\phi\circ\mu_c$ holds, it follows from Propositions 3.2 and 3.3 that 
${\rm hol}_c:(\mathcal A_P^{H^s},g_s)\to(G,g_G)$ is a homothetic submersion of coefficient $a$.  
We consider the case of $s=0$.  Then it is clear that $\mu_c$ is an affine map and the operator norm of the linear part is equal to $a$ 
(and hence it is bounded).  Therefore, the fibres of $\mu_c$ are 
affine subspaces of $\mathcal A_P^{H^0}$ and hence they are totally geodesic in $(\mathcal A_P^{H^0},g_0)$.  
On the other hand, according to Proposition 3.6, $\phi:(H^0([0,1],\mathfrak g),g_0^o)\to(G,g_G)$ has a Riemannian submersion with minimal regularizable fibres.  
Therefore, ${\rm hol}_c:(\mathcal A_P^{H^0},g_0)\to (G,g_G)$ is a homothetic submersion of coefficient $a$ with minimal regularizable fibres.  
Thus the statement (i) of Theorem A is shown.  

Let $M$ be an equifocal submanifold in $(G,g_G)$.  Since ${\rm hol}_c:(\mathcal A_P^{H^0},g_0)\to (G,\frac{1}{a^2}g_G)$ is a Riemannian submersion with minimal regularizable 
fibres and $M$ is equifocal in $(G,\frac{1}{a^2}g_G)$, it follows from Theorem 3.7 that ${\rm hol}_c^{-1}(M)$ is isoparametric in $(\mathcal A_P^{H^0},g_0)$.  
The $(\mathcal G_P^{H^1})_{c(0)}$-invariance of ${\rm hol}_c^{-1}(M)$ follows from $(2.8)$.  
Also, in the case where $M$ is ${\rm Ad}(G)$-invariant, the $\mathcal G_P^{H^1}$-invariance of ${\rm hol}_c^{-1}(M)$ follows from $(2.9)$.  \qed

\end{document}